\documentclass[12pt,reqno]{amsart}
\usepackage{amsmath, amsfonts, amssymb, amsthm, hyperref}
\usepackage{bm}
\allowdisplaybreaks[4]
\def\l{\left}
\def\r{\right}
\def\bg{\bigg}
\def\({\bg(}
\def\){\bg)}
\def\t{\text}
\def\f{\frac}

\def\ch{{\rm ch}}

\def\ls{\leqslant}
\def\gs{\geqslant}
\def\se {\subseteq}
\def\sm{\setminus}
\def\bi{\binom}

\def\eq{\equiv}

\def\da{\delta}

\def\Proof{\noindent{\it Proof}}

\def\Z{\mathbb Z}
\def\C{\mathbb C}
\def\N{\mathbb N}
\def\Q{\mathbb Q}

\def\1{{\bf 1}}

\def\pmod #1{\ ({\rm{mod}}\ #1)}

\def\<{\langle}
\def\>{\rangle}
\theoremstyle{plain}
\newtheorem{theorem}{Theorem}[section]
\newtheorem{lemma}{Lemma}
\newtheorem{corollary}{Corollary}
\newtheorem{conjecture}{Conjecture}

\theoremstyle{definition}

\theoremstyle{remark}

\makeatletter
\@namedef{subjclassname@2020}{%
  \textup{2020} Mathematics Subject Classification}
\makeatother
 \vspace{4mm}

\begin{document}
\hbox{Preprint}
\medskip

\title[On two new kinds of restricted sumsets]
{On two new kinds of restricted sumsets}
\author{Han Wang}
\address {(Han Wang) Department of Mathematics, Nanjing
University, Nanjing 210093, People's Republic of China}
\email{hWang@smail.nju.edu.cn}

\author{Zhi-Wei Sun}
\address {(Zhi-Wei Sun, corresponding author) Department of Mathematics, Nanjing
University, Nanjing 210093, People's Republic of China}
\email{zwsun@nju.edu.cn}

\keywords{Additive combinatorics, sumset, field, the polynomial method, torsion-free abelian group.
\newline \indent Supported by the Natural Science Foundation of China (grant no. 11971222).}
\begin{abstract}
Let $A_1,\ldots,A_n$ be finite subsets of an additive abelian group $G$ with $|A_1|=\cdots=|A_n|\gs2$.
Concerning the two new kinds of restricted sumsets
$$L(A_1,\ldots,A_n)=\{a_1+\cdots+a_n:\ a_1\in A_1,\ldots,a_n\in A_n,\ \text{and}\ a_i\not=a_{i+1}
\ \text{for}\ 1\ls i<n\}$$
and
$$C(A_1,\ldots,A_n)=\{a_1+\cdots+a_n:\ a_i\in A_i\ (1\ls i\ls n),\ \text{and}\ a_i\not=a_{i+1}
\ \text{for}\ 1\ls i<n,\ \text{and}\ a_n\not=a_1\}$$
recently introduced by the second author, when $G$ is the additive group of a field we obtain lower bounds for $|L(A_1,\ldots,A_n)|$
and $|C(A_1,\ldots,A_n)|$ via the polynomial method. Moreover,
when $G$ is torsion-free and $A_1=\cdots=A_n$, we determine completely when $|L(A_1,\ldots,A_n)|$ or
$|C(A_1,\ldots,A_n)|$ attains its lower bound.
\end{abstract}
\maketitle

\section{Introduction}
\setcounter{lemma}{0}
\setcounter{theorem}{0}
\setcounter{equation}{0}
\setcounter{conjecture}{0}
\setcounter{remark}{0}
\setcounter{corollary}{0}

Let $G$ be an additive group. For finite subsets $A_1,\ldots,A_n$ of $G$, their sumset
is defined by
$$A_1+\cdots+A_n=\{a_1+\cdots+a_n:\ a_1\in A_1,\ldots,a_n\in A_n\},$$
and their sumset with distinct summands is given by
$$A_1\dotplus\cdots\dotplus A_n=\{a_1+\cdots+a_n:\ a_1\in A_1,\ldots,a_n\in A_n,
\ \t{and}\ a_1,\ldots,a_n\ \t{are pairwise distinct}\}.$$
When $A_1=\cdots=A_n=A$, we use $nA$ to mean $A_1+\cdots+A_n$, and $n^\wedge A$ to stand for $A_1\dotplus
\cdots\dotplus A_n$.

The famous Erd\H os-Heilbronn conjecture \cite{EH} posed in 1964 asserts that
for any prime $p$ and $A\se\Z/p\Z$, we have
$$|2^\wedge A|\gs\min\{p,2|A|-3\}.$$
This was first confirmed by Dias da Silva and Hamidoune \cite{DH} in 1994 via exterior algebra.
In 2009 P. Balister and J. P. Wheeler \cite{BW} obtained the following extension of the Erd\H os-Heilbronn conjecture
to any finite group $G$ with $|G|>1$: For any finite subset $A$ of $G$ we have
$$|2^\wedge A|\gs\min\{p(G),2|A|-3\},$$
where $p(G)$ is the minimum of orders of nonzero elements of $G$ (which is $+\infty$ if each nonzero element of $G$ has infinite order).

For a subset $A$ of an additive group $G$, note that
$$2^\wedge A=\{a_1+a_2:\ a_1,a_2\in A\ \t{and}\ a_1\not=a_2\}$$
and
$$3^\wedge A=\{a_1+a_2+a_3:\ a_1,a_2,a_3\in A\ \t{and}\ a_1\not=a_2\not=a_3\not=a_1\}.$$
Motivated by this, and linear and circular permutations,
recently the second author \cite{OEIS} introduced two new kinds of sumsets. Namely,
for finite subsets $A_1,\ldots,A_n$ of an additive group $G$, Sun defined
$$L(A_1,\ldots,A_n)=\{a_1+\cdots+a_n:\ a_i\in A_i\ \t{for}\ i=1,\ldots,n,
\ \t{and}\ a_{i}\not=a_{i+1}\ \t{for all}\ 0<i<n\}$$
and
$$C(A_1,\ldots,A_n)=\{a_1+\cdots+a_n:\ a_i\in A_i\ \t{for}\ i=1,\ldots,n,
\ a_{i}\not=a_{i+1}\ \t{for all}\ 0<i<n,\ \t{and}\ a_n\not= a_1\}.$$
 Clearly
$$L(A_1,A_2)=C(A_1,A_2)=A_1\dotplus A_2,\ \t{and}\ C(A_1,A_2,A_3)=A_1\dotplus A_2\dotplus A_3.$$
When $A_1=\cdots=A_n=A$, we simply write $n\,\tilde{}A$ to denote $L(A_1,\ldots,A_n)$,
and $n^\circ A$ to denote $C(A_1,\ldots,A_n)$.

Recently, the second author \cite{OEIS} made the following general conjecture.

\begin{conjecture}\label{ConjG} Let $G$ be an additive group with $|G|>1$, and let $A_1,\ldots,A_n\ (n>1)$
be finite subsets of $G$ with $|A_i|>1$ for all $i=1,\ldots,n$. Then
\begin{equation}\label{L}
|L(A_1,\ldots,A_n)|\gs\min\l\{p(G),\, |A_1|+\cdots+|A_n|-2n+1+\{n\}_2\r\}
\end{equation}
and
\begin{equation}\label{C}
|C(A_1,\ldots,A_n)|\gs\min\{p(G),\, |A_1|+\cdots+|A_n|-2n+(-1)^n(1+\{n\}_2)\},
\end{equation}
where $\{n\}_2$ denotes the least nonnegative residue of $n$ modulo $2$.
\end{conjecture}

This conjecture is motivated by Sun's following observation: For any integer $n>1$ and $A=\{0, \ldots,k-1\}$ with $k\in\mathbb{Z}^+=\{1,2,3,\ldots\}$, we have
\[
|n\,\tilde{} A|=n|A|-2n+1+\{n\}_2
\]
and
\[
|n^\circ A|=n|A|-2n+(-1)^n(1+\{n\}_2).
\]
By Sun \cite[Corollary 1.5]{S01}, if $A_i\se\Z$ and $|A_i|\gs3$ for all $i=1,\ldots,n$,
then $|C(A_1,\ldots,A_n)|\gs\sum_{i=1}^n|A_i|-3n+1$.

For a field $F$, we let $\ch(F)$ be the characteristic of $F$.
Clearly,
$$p(F)=\begin{cases}p&\t{if}\ \ch(F)\ \t{is a prime}\ p,
\\+\infty&\t{if}\ \ch(F)=0.\end{cases}$$
Applying Theorem 1.3 of Sun and Zhao \cite{SunZhao} with
$$P(x_1,\ldots,x_n)=(x_1-x_2)(x_2-x_3)\cdots(x_{n-1}-x_n)(x_n-x_1)$$
or
$$P(x_1,\ldots,x_n)=(x_1-x_2)(x_2-x_3)\cdots(x_{n-1}-x_n),$$
we see that if $A_1,\ldots,A_n$ are finite subsets of a field $F$
with $|A_i|>2$ for all $i=1,\ldots,n$ then
$$|C(A_1,\ldots,A_n)|\gs\min\{p(F)-n,\,|A_1|+\cdots+|A_n|-3n+1\}$$
and
$$|L(A_1,\ldots,A_n)|\gs\min\{p(F)-n+1,\,|A_1|+\cdots+|A_n|-3n+3\}.$$

 If $A_1,A_2,A_3$ are finite nonempty subsets of a field $F$,  then in the spirit of \cite[Theorem 3.2]{ANR}, the cardinality of $C(A_1,A_2,A_3)=A_1\dotplus A_2\dotplus A_3$ is
at least
$$\min\l\{p(F),\,1+\sum_{i=1}^3(|A_i|-3)\r\}=\min\{p(F),\,|A_1|+|A_2|+|A_3|-2\times3+(-1)^3(1+\{3\}_2)\},$$
which is essentially the inequality \eqref{C} with $n=3$ and $G=F$.

When $G$ is the additive group of a field,
our following theorem confirms the inequality \eqref{C} in the case $n=3$ and $|A_1|=|A_2|=|A_3|$.

\begin{theorem}\label{Th1.1} {\rm (i)} Let $F$ be any field, and let $A_1,A_2,A_3$ be finite subsets of $F$ with $|A_1|,|A_3|\gs2$ and $|A_2|-|A_1|\in\{0,1\}$.
Then we have
\begin{equation}\label{L3}|L(A_1,A_2,A_3)|\gs\min\{p(F),|A_1|+|A_2|+|A_3|-4\}.
\end{equation}

{\rm (ii)} Let $G$ be a torsion-free additive abelian group,
 and let $A_1,A_2,A_3$ be finite subsets of $G$ with cardinality $2\ls |A_1|\ls|A_2|\ls |A_3|$.
Then we have
\begin{equation}\label{L3G}|L(A_1,A_2,A_3)|\gs|A_1|+|A_2|+|A_3|-4.
\end{equation}
\end{theorem}

\begin{corollary} \label{Cor1} Let $p$ be any prime. For any $A\se\mathbb F_p$ with $|A|\gs\lfloor p/3\rfloor+2$,
each element of $\mathbb F_p$ can be written as $a_1+a_2+a_3$ with $a_1,a_2,a_3\in A$ and $a_1\not=a_2\not=a_3$.
\end{corollary}

 Our next two theorems deal with Conjecture \ref{ConjG} when $G$ is the additive group of a field,
 and $|A_1|=\cdots=|A_n|$.

\begin{theorem}\label{Th1.2} Let $n$ be any positive even integer,
and let $A_1,\dots,A_n$ be subsets of a field $F$ with $|A_1|=\cdots=|A_n|\gs2$.
Suppose that $p(F)>\sum_{i=1}^n|A_i|-2n$. Then
\begin{equation}\label{C-even}|L(A_1,\ldots,A_n)|\gs |C(A_1,\ldots,A_n)|\gs\sum_{i=1}^n|A_i|-2n+1.
\end{equation}
\end{theorem}

\begin{theorem}\label{Th1.3} Let $n$ be any positive odd integer,
and let $A_1,\dots,A_n$ be subsets of a field $F$ with $|A_1|=\cdots=|A_n|\gs2$.
Suppose that $p(F)>\sum_{i=1}^n|A_i|-2n+1$. Then
\begin{equation}\label{L-odd}|L(A_1,\ldots,A_n)|\gs\sum_{i=1}^n|A_i|-2n+2.
\end{equation}
\end{theorem}

We will prove Theorem \ref{Th1.1} and Corollary \ref{Cor1},
and Theorem \ref{Th1.2}--\ref{Th1.3} in Sections 2 and 3 respectively, via the so-called polynomial method involving the famous Combinatorial Nullstellensatz of N. Alon \cite{CN}.

In 1995, M. B. Nathanson \cite{N95} proved that for any finite subset $A$ of $\Z$ with $|A|\gs5$ and
$n\in\{2,\ldots,|A|-2\}$ we have $|n^\wedge A|\gs n|A|-n^2+1$, and equality holds if and only if $A$
is an AP (arithmetic progression).

For any torsion-free abelian group $G$,  we confirm Conjecture \ref{ConjG}
in the case $A_1=\cdots=A_n$. Namely, we have the following theorem.

\begin{theorem}\label{Th1} Let $A$ be any finite subset of
a torsion-free additive abelian group $G$ with $|A|\gs2$. For any integer $n>1$, we have
\begin{equation}\label{-}
|n\,\tilde{} A|\gs n|A|-2n+1+\{n\}_2
\end{equation}
and
\begin{equation}\label {o}
|n^\circ A|\gs n|A|-2n+(-1)^n(1+\{n\}_2).
\end{equation}
\end{theorem}

When $G$ is a torsion-free abelian group, we determine completely when equality in \eqref{-} or \eqref{o} holds.

\begin{theorem}\label{Th2} Let $n>2$ be an integer, and let $A$ be any finite subset of
an additive torsion-free abelian group $G$ with $|A|\gs3$. Then
\begin{equation}\label{-eq}
|n\,\tilde{} A|=n|A|-2n+1+\{n\}_2
\end{equation}
if and only $A$ is an AP (arithmetic progression).
Also,
\begin{equation}\label{oeq}
|n^\circ A|=n|A|-2n+(-1)^n(1+\{n\}_2)
\end{equation}
if and only if $A$ is an AP,  or $k=3$ and $n=5$.
\end{theorem}

\section{Proofs of Theorem \ref{Th1.1} and Corollary \ref{Cor1}}
\setcounter{lemma}{0}
\setcounter{theorem}{0}
\setcounter{equation}{0}
\setcounter{conjecture}{0}
\setcounter{remark}{0}
\setcounter{corollary}{0}

For a polynomial $P(x_1,\ldots,x_n)$ over a field, and nonnegative integers $k_1,\ldots,k_n$
as usual we write $[x^{k_1}\cdots x_n^{k_n}]P(x_1,\ldots,x_n)$
to mean the coefficient of $x^{k_1}\cdots x_n^{k_n}$ in $P(x_1,\ldots,x_n)$.

The following lemma is essentially Theorem 2.1 of \cite{ANR} although its original form deals with $F=\Z/p\Z$ with $p$ prime, the reader may also consult \cite{S07} for a proof
via Alon's Combinatorial Nullstellensatz \cite{CN}.

\begin{lemma}\label{Lem1} Let $F$ be any field, and let $A_1,\ldots, A_n$ be finite nonempty subsets of $F$. Let $f(x_1,\ldots,x_n)\in F[x_1,\ldots,x_n]\sm\{0\}$ with $\deg{f}\ls\sum_{i=1}^n(|A_i|-1)$. If
\[
[x_1^{|A_1|-1}\ldots x_n^{|A_n|-1}]f(x_1,\ldots,x_n)(x_1+\cdots+x_n)^{\sum_{i=1}^n(|A_i|-1)-\deg f}\ne 0,
\]
then we have
\[
|\{a_1+\cdots+a_n:a_i\in A_i,f(a_1,\ldots,a_n)\ne 0\}|\gs\sum_{i=1}^n(|A_i|-1)-\deg f+1.\]
\end{lemma}

The following lemma is essentially due to N. Alon, M.B. Nathanson and I.Z. Ruzsa \cite{ANR95,ANR}.

\begin{lemma}\label{AB} If $A$ and $B$ are finite subsets of a field $F$ with $0<|A|<|B|$, then
\begin{equation}|A\dotplus B|\gs\min\{p(F),\,|A|+|B|-2\}.
\end{equation}
\end{lemma}

\begin{lemma} For any $k_1,k_2,k_3\in\Z^+$, we have
\begin{equation}\label{coeff}\begin{aligned}
&[x_1^{k_1}x_2^{k_2}x_3^{k_3}](x_1-x_2)(x_2-x_3)(x_1+x_2+x_3)^{k_1+k_2+k_3-2}
\\=&\f{(k_1+k_2+k_3-2)!}{k_1!k_2!k_3!}(k_2+(k_2-k_1)(k_3-k_2)).
\end{aligned}\end{equation}
\end{lemma}
\Proof. Note that $(x_1-x_2)(x_2-x_3)=x_1x_2-x_1x_3+x_2x_3-x_2^2$.
With the aid of the multi-nomial theorem, we see that
\begin{align*}&[x_1^{k_1}x_2^{k_2}x_3^{k_3}](x_1x_2-x_1x_3+x_2x_3-x_2^2)(x_1+x_2+x_3)^{k_1+k_2+k_3-2}
\\=&\ \bi{k_1+k_2+k_3-2}{k_1-1,k_2-1,k_3}-\bi{k_1+k_2+k_3-2}{k_1-1,k_2,k_3-1}
+\bi{k_1+k_2+k_3-2}{k_1,k_2-1,k_3-1}-\bi{k_1+k_2+k_3-2}{k_1,k_2-2,k_3}
\\=&\ \f{(k_1+k_2+k_3-2)!}{k_1!k_2!k_3!}(k_1k_2-k_1k_3+k_2k_3-k_2(k_2-1))
\\=&\ \f{(k_1+k_2+k_3-2)!}{k_1!k_2!k_3!}(k_2+(k_2-k_1)(k_3-k_2)).
\end{align*}
This proves \eqref{coeff}. \qed

\begin{lemma}\label{Lem2} Let $F$ be any field, and let $A_1,A_2,A_3$ be
finite subsets of $F$ with $|A_1|,|A_3|\gs2$ and $|A_2|-|A_1|\in\{0,1\}$.
Suppose that $p(F)\gs\sum_{i=1}^3|A_i|-4$.
Then we have
\begin{equation}|L(A_1,A_2,A_3)|\gs\sum_{i=1}^3|A_i|-4.
\end{equation}
\end{lemma}
\Proof.
Set $k_i=|A_i|-1$ for $i=1,2,3$. In light of Lemma \ref{Lem1}, it suffices to prove that
$he\not=0$, where $e$ is the identity of $F$, and $h$ is the coefficient
of $x_1^{k_1}x_2^{k_2}x_3^{k_3}$ in the polynomial
$$(x_1-x_2)(x_2-x_3)(x_1+x_2+x_3)^{k_1+k_2+k_3-2}\in\Z[x].$$

Let $\da=|A_2|-|A_1|$. Then $k_2-k_1=\da\in\{0,1\}$. In light of \eqref{coeff}, we have
$$h=\f{(k_1+k_2+k_3-2)!k_{2+\da}}{k_1!k_2!k_3!}.$$
Since $p(F)\gs \sum_{i=1}^3(|A_i|-1)-1=k_1+k_2+k_3-1>k_{2+\da}$, we clearly have $he\not=0$ as desired.
This concludes the proof. \qed

\medskip
\noindent{\it Proof of Theorem \ref{Th1.1}}.
(i) When $p(F)\gs \sum_{i=1}^3|A_i|-4$, the conclusion follows from Lemma \ref{Lem2}.

Below we assume that $p(F)=p<\sum_{i=1}^3|A_i|-4$.
\medskip

{\it Case} 1. $p=2$.

In this case, $|A_1|+|A_2|>p+4-|A_3|=6-|A_3|$.

Suppose $|A_3|=2$. Then $|A_1|+|A_2|>4$ and hence $|A_2|\gs3$. Let $a_1$ and $a_1'$ be two distinct
elements of $A_1$, and choose $a_2\in A_2$ different from $a_1$ and $a_1'$. Also, take $a_3\in A_3$
with $a_3\not=a_2$. Then $a_1+a_2+a_3$ and $a_1'+a_2+a_3$ are distinct elements of $L(A_1,A_2,A_3)$
and hence
$$|L(A_1,A_2,A_3)|\gs2=p=\min\{p,\,|A_1|+A_2|+|A_3|-4\}.$$

Now assume $|A_3|\gs3$. Take $a_1\in A_1$ and choose $a_2\in A_2$ with $a_2\not=a_1$. As $|A_3|\gs3$, there are two distinct elements $a_3$ and $a_3'$ of $A_3$
different from $a_2$. Thus $a_1+a_2+a_3$ and $a_1+a_2+a_3'$ are distinct elements of $L(A_1,A_2,A_3)$,
and hence $|L(A_1,A_2,A_3)|\gs 2=\min\{p,|A_1|+|A_2|+|A_3|-4\}$.
\medskip

{\it Case} 2. $p\not=2$ and $|A_3|=2$.

In this case, $|A_1|+|A_2|>p+4-|A_3|=p+2$.
Let $b$ be an element of $A_3$. We claim that
\begin{equation}\label{claim}
|A_1\dotplus (A_2\sm\{b\})|\gs p.
\end{equation}

{\it Subcase} 2.1. $b\in A_2$.
\medskip

By Lemma \ref{AB}, if $|A_1|=|A_2|$ then
$$|A_1\dotplus (A_2\sm\{b\})|\gs\min\{p,|A_1|+|A_2\sm\{b\}|-2\}=\min\{p,|A_1|+|A_2|-3\}=p.$$
When $|A_2|=|A_1|+1$, we take $a\in A_1$ and then by Lemma \ref{AB} we get
\begin{align*}|A_1\dotplus A_2\sm \{b\}|
&\gs |(A_1\sm\{a\})\dotplus (A_2\sm\{b\})|
\\&\gs\min\{p,\, |A_1\sm\{a\}|+|A_2\sm\{b\}|-2\}=\min\{p,\, |A_1|+|A_2|-4\}.
\end{align*}
As $|A_1|+|A_2|-4\gs p-1$ and $|A_2|\not\eq |A_1|\pmod2$, we must have $|A_1|+|A_2|\gs p+4$
and hence $|A_1\dotplus A_2\sm \{b\}|\gs p$.
\medskip

{\it Subcase} 2.2. $b\not\in A_2$.
\medskip

When $|A_1|=|A_2|$, we take $c\in A_2$, and hence by Lemma \ref{AB} we have
$$|A_1\dotplus (A_2\sm\{b\})|\gs |A_1\dotplus (A_2\sm\{c\})|\gs\min\{p,|A_1|+|A_2\sm\{c\}|-2\}
=\min\{p,|A_1|+|A_2|-3\}=p.$$
If $|A_2|=|A_1|+1$, then
by Lemma \ref{AB} we have
$$|A_1\dotplus (A_2\sm\{b\})|\gs |A_1\dotplus A_2|\gs\min\{p,|A_1|+|A_2|-2\}=p.$$

By our discussion of subcases 2.1 and 2.2, we see that \eqref{claim} holds.
For any $x\in A_1\dotplus (A_2\sm\{b\})$, we may write $x=a_1+a_2$ with $a_1\in A_1$, $a_2\in A_2$
and $a_1\not=a_2\not= b$ and thus $x+b\in L(A_1,A_2,A_3)$.
Combining this with \eqref{claim}, we obtain
$$|L(A_1,A_2,A_3)|\gs  |A_1\dotplus (A_2\sm\{b\})|\gs p=\min\{p,\,|A_1|+|A_2|+|A_3|-4\}.$$

{\it Case} 3. $p\gs3$ and $|A_3|\gs3$.

If $|A_2|\ls (p+1)/2$, then $3\ls p+4-(|A_1|+|A_2|)<|A_3|$ and so we may take $A_3'\se A_3$
with $|A_3'|=p+4-(|A_1|+|A_2|)$, hence by Lemma \ref{Lem2} we have
$$|L(A_1,A_2,A_3)|\gs|L(A_1,A_2,A_3')|\gs |A_1|+|A_2|+|A_3'|-4=p=\min\{p,\,|A_1|+|A_2|+|A_3|-4\}.$$

Below we assume that  $|A_2|>(p+1)/2$.   We  may take $A_1'\se A_1$ and $A_2'\se A_2$
with $|A_1'|=|A_2'|=(p+1)/2$, and also take $A_3'\se A_3$ with $|A_3'|=3$.
Note that $|A_1'|+|A_2'|+|A_3'|-4=p$. Applying Lemma \ref{Lem2} we obtain
$$|L(A_1,A_2,A_3)|\gs|L(A_1',A_2',A_3')|\gs p=\min\{p,|A_1|+|A_2|+|A_3|-4\}.$$

(ii) Now we turn to prove part (ii) of Theorem \ref{Th1.1}. Note that the subgroup of $G$
generated by $A_1\cup A_2\cup A_3$ is a finitely generated torsion-free abelian group
which is isomorphic to $\Z^r$ for some positive integer $r$.
Let $K$ be algebraic number field with $[K:\Q]=r$, and let $O_K$ be the ring of all algebraic integers.
It is known that $O_K\cong \Z^r$ (cf. \cite[p.\,175]{IR}). Thus, without loss of generality, we may simply let $G$
be the additive group of the field $\C$.

Let $k_i=|A_i|-1$ for $i=1,2,3$. Then $1\ls k_1\ls k_2\ls k_3$ and hence
$$[x_1^{k_1}x_2^{k_2}x_3^{k_3}](x_1-x_2)(x_2-x_3)(x_1+x_2+x_3)^{k_1+k_2+k_3-2}>0$$
by \eqref{coeff}. Thus, applying Lemma \ref{Lem1} we see that
$$|L(A_1,A_2,A_3)|\gs \sum_{i=1}^3 k_i-1=|A_1|+|A_2|+|A_3|-4.$$

In view of the above, we have completed the proof of Theorem \ref{Th1.1}. \qed
\medskip

\noindent{\it Proof of Corollary} \ref{Cor1}. Applying Theorem \ref{Th1.1}(i) with $F=\mathbb F_p$
and $A_1=A_2=A_3=A$, we get
$$|3\,\tilde{} A|\gs\min\{p,3|A|-4\}.$$
As $|A|\gs\lfloor p/3\rfloor +2$, we have
$3|A|-4\gs3\lfloor p/3\rfloor+2\gs p$. Thus $|3\,\tilde{} A|\gs p$ and hence $3\,\tilde{} A=\mathbb F_p$.
This concludes the proof. \qed

\section{Proofs of Theorems \ref{Th1.2} and \ref{Th1.3}}
\setcounter{lemma}{0}
\setcounter{theorem}{0}
\setcounter{equation}{0}
\setcounter{conjecture}{0}
\setcounter{remark}{0}
\setcounter{corollary}{0}

The following lemma is essentially due to Q.-H. Hou and Z.-W. Sun \cite{HS}, and
a generalization was given by Sun and Y.-N. Yeh \cite[Lemma 2.1]{SY}.

\begin{lemma}\label{P*}
 Let  $$P(x_1,\ldots,x_n)=\sum_{j_1,\ldots,j_n\gs0\atop j_1+\cdots+j_n=m} c_{j_1,\ldots,j_n}
 x_1^{j_1}\cdots x_n^{j_n}\in\C[x_1,\ldots,x_n]$$
 and
 $${\mathcal L}(P)(x)=\sum_{j_1,\ldots,j_n\gs0\atop j_1+\cdots+j_n=m} c_{j_1,\ldots,j_n}
 (x)_{j_1}\cdots (x)_{j_n}.$$
 Suppose that $0\ls\deg P\ls kn$
 with $k\in\N$. Then
 $$[x_1^{k}\cdots x_n^{k}]P(x_1,\ldots,x_n)(x_1+\cdots+x_n)^{kn-\deg P}
=\f{(kn-\deg P)!}{(k!)^n}{\mathcal L}(P)(k).$$
\end{lemma}

\begin{lemma}\label{Lem2.5} Let
$$P_n(x_1,\ldots,x_n)=(x_1-x_2)(x_2-x_3)\cdots(x_{n-1}-x_n)(x_n-x_1)$$
with $n$ even. Then we have
\begin{equation}{\mathcal L}(P_n)(x)=2x^{n/2}.
\end{equation}
\end{lemma}
\Proof.
Write
$$P_n(x_1,\ldots,x_n)=\sum_{I,J\se\{1,\ldots,n\}\atop I\cap J=\emptyset\ \&\ |I|=|J|}
c(I,J)\prod_{i\in I} x_i^2\times\prod_{j\in J}x_j^0\times\prod_{k=1\atop k\not\in I\cup J}^nx_k$$
with $c(I,J)\in\Z$. Then
\begin{equation}
\label{LP} {\mathcal L}(P_n)(x)=\sum_{I,J\se\{1,\ldots,n\}\atop I\cap J=\emptyset\ \&\ |I|=|J|}
\prod_{i\in I}(x)_{2}\times\prod_{k=1\atop k\not\in I\cup J}^n(x)_1
=\sum_{m=0}^{n/2}c_m(x)_2^m(x)_1^{n-2m}=\sum_{m=0}^{n/2}c_mx^{n-m}(x-1)^m,
\end{equation}
where
$$c_m=\sum_{I,J\se\{1,\ldots,n\}\atop I\cap J=\emptyset\ \&\ |I|=|J|=m} c(I,J).$$

Let $I,J\se\{1,\ldots,n\}$ with $I\cap J=\emptyset$ and $|I|=|J|=m$.  Suppose that  $j,k\in J$, $j<k$ and $s\not\in J$ for all $j<s<k$. Then there are only $k-j-1$ choices of the corresponding terms chosen from $$(x_j-x_{j+1})(x_{j+1}-x_{j+2})\ldots(x_{k-2}-x_{k-1})(x_{k-1}-x_k),$$ which are
$$\prod_{j<r\ls s}(-x_r)\times\prod_{s\ls t<k}x_t=(-1)^{s-j}x_s^2\prod_{j<r<k\atop r\not=s}x_r
\ \ (j<s<k).$$
Note that
$$\sum_{j<s<k}{\mathcal L}\((-1)^{s-j}x_s^2\prod_{j<r<k\atop r\not=s}x_r\)
=\sum_{j<s<k}(-1)^{s-j}x(x-1)\prod_{j<r<k\atop r\not=s}x=(\{k-j\}_2-1)(x-1)x^{k-j-1}$$
which vanishes if $j\not\eq k\pmod 2$.

Let $J\se\{1,\ldots,n\}$ with $|J|=m$. Write $J=\{j_1,j_2,\ldots,j_m\}$ with $j_1<\cdots<j_m$.
In the spirit of the last paragraph, $j_1\eq j_2\eq\cdots\eq j_m\pmod 2$, and for the polynomial
$P_J$ given by
\begin{align*}&\prod_{i=1}^{m-1}\(\sum_{j_i<s<j_{i+1}}\prod_{j_i<r\ls s}(-x_r)\times\prod_{s\ls t<j_{i+1}} x_t\)
\\&\times \(\sum_{j_m<s\ls n}\prod_{j_m<r\ls s}(-x_r)\times\prod_{s\ls t\ls n\atop \t{or}\ 1\ls t<j_1}x_t
+\sum_{1\ls s<j_1}\prod_{j_m<r\ls n\atop\t{or}\ 1\ls r\ls s}(-x_r)\times\prod_{s\ls t<j_1}x_t\)
\end{align*}
we have
$${\mathcal L}(P_J)(x)=\prod_{i=1}^{m-1}\l(-(x-1)x^{j_{i+1}-j_i-1}\r)
\times\l(-(x-1)x^{(n+j_1)-j_m-1}\r)=(-1)^m(x-1)^mx^{n-m}.$$
Thus
$$\sum_{I\se\{1,\ldots,n\}\sm J\atop |I|=m}c(I,J)=(-1)^m.$$

By the last paragraph, we have
\begin{align*}c_m&=\sum_{J\se\{2s:\ s=1,\ldots,n/2\}\atop |J|=m}\sum_{I\se\{1,\ldots,n\}\sm J\atop |I|=m}c(I,J)
+\sum_{J\se\{2s-1:\ s=1,\ldots,n/2\}\atop |J|=m}\sum_{I\se\{1,\ldots,n\}\sm J\atop |I|=m}c(I,J)
\\&=\bi{n/2}m(-1)^m+\bi{n/2}m(-1)^m = (-1)^m 2\bi{n/2}m.
\end{align*}
Combining this with \eqref{LP} we get
$$
{\mathcal L}(P_n)(x)=\sum_{m=0}^{n/2}(-1)^m 2\bi{n/2}mx^{n-m}(x-1)^m=2x^{n/2}$$
by the binomial theorem. This concludes our proof. \qed

\medskip
\noindent{\it Proof of Theorem} \ref{Th1.2}. It is apparent that
$$L(A_1,\ldots,A_n)\supseteq C(A_1,\ldots,A_n).$$
So, we only to show the second inequality of \eqref{C-even}.

Let $k=|A_1|-1=\cdots=|A_n|-1$.
If $k=1$, then the second equality of \eqref{C-even} holds trivially.

Below we assume that $k\gs2$.
By Lemmas \ref{P*} and \ref{Lem2.5}, we have the identity
\begin{equation}\label{even-n}
[x_1^{k}\ldots x_n^{k}](x_1-x_2)\cdots(x_{n-1}-x_n)(x_n-x_1)(x_1+\cdots+x_n)^{(k-1)n}=\f{((k-1)n)!}{k!^n}2k^{n/2}.
\end{equation}
Let $h$ denote the right-hand side of \eqref{even-n} which is an integer. Let $e$ be the identity of $F$. As $p(F)>(k-1)n\gs2$, we see that the coefficient of $x_1^{k}\ldots x_n^{k}$
in the polynomial
$$(x_1-x_2)\cdots(x_{n-1}-x_n)(x_n-x_1)(x_1+\cdots+x_n)^{(k-1)n}\in F[x_1,\ldots,x_n]$$
coincides with $he$ which is nonzero. Therefore, by Lemma \ref{Lem1} we have
$$|C(A_1,\ldots,A_n)|\gs(k-1)n+1=\sum_{i=1}^n|A_i|-2n+1.$$
This concludes our proof. \qed

\begin{lemma}\label{Lem-Q} For any odd integer $n>1$ and the polynomial
$$Q_n(x_1,\ldots,x_n)=(x_1-x_2)(x_2-x_3)\cdots(x_{n-1}-x_n),$$
we have
\begin{equation}\label{LQ} {\mathcal L}(Q_n)(x)=x^{(n-1)/2}.
\end{equation}
\end{lemma}
\Proof. We use induction on $n$.

Clearly
$$Q_3(x_1,x_2,x_3)=(x_1-x_2)(x_2-x_3)=x_1x_2-x_1x_3+x_2x_3-x_2^2$$
and hence
$${\mathcal L}(Q_3)(x)=xx-xx+xx-(x)_2=x^2-x(x-1)=x^{(3-1)/2}.$$
Thus \eqref{LQ} holds when $n=3$.

Now let $n$ be an odd integer greater than $3$, and assume that
${\mathcal L}(Q_{n-2})(x)=x^{(n-3)/2}$. Write
$$Q_n(x_1,\ldots,x_n)=\sum_{i_1,\ldots,i_n\in\N}c_{i_1,\ldots,i_n}x_1^{i_1}\cdots x_n^{i_n}
\ \t{with}\ c_{i_1,\ldots,i_n}\in\Z,$$
and define
$$P_{n+1}(x_1,\ldots,x_n,x_{n+1})=Q_n(x_1,\ldots,x_n)(x_n-x_{n+1})(x_{n+1}-x_1).$$
Then
\begin{align*}P_{n+1}(x_1,\ldots,x_{n+1})&=\sum_{i_1,\ldots,i_n\in\N}
c_{i_1,\ldots,i_n}x_1^{i_1}\cdots x_n^{i_n}(x_nx_{n+1}-x_1x_n+x_1x_{n+1}-x_{n+1}^2)
\\&=\sum_{i_1,\ldots,i_n\in\N}
c_{i_1,\ldots,i_n}\l(x_1^{i_1}\cdots x_{n-1}^{i_{n-1}}x_n^{i_n+1}x_{n+1}
-x_1^{i_1+1}x_2^{i_2}\cdots x_{n-1}^{i_{n-1}}x_n^{i_n+1}\r)
\\&\ \ \ \ +\sum_{i_1,\ldots,i_n\in\N}
c_{i_1,\ldots,i_n}\l(x_1^{i_1+1}x_2^{i_2}\cdots x_n^{i_n}x_{n+1}
-x_1^{i_1}\cdots x_n^{i_n}x_{n+1}^2\r)
\end{align*}
and hence
$${\mathcal L}(P_{n+1})(x)=\sum_{i_1,\ldots,i_n\in\N}c_{i_1,\ldots,i_n}R(i_1,\ldots,i_n,x),$$
where
\begin{align*}R(i_1,\ldots,i_n,x)&=(x)_{i_1}\cdots(x)_{i_{n-1}}(x)_{i_n+1}x
-(x)_{i_1+1}(x)_{i_2}\cdots(x)_{i_{n-1}}(x)_{i_n+1}
\\&\ \ \ \ +(x)_{i_1+1}(x)_{i_2}\cdots (x)_{i_n}x-(x)_{i_1}\cdots(x)_{i_n}(x)_2
\\&= (x)_{i_1}\cdots(x)_{i_n}\l(x(x-i_n)-(x-i_1)(x-i_n)+x(x-i_1)-x(x-1)\r)
\\&=(x)_{i_1}\cdots(x)_{i_n}(x-i_1i_n).
\end{align*}
Thus
\begin{equation}
\label{Relation}{\mathcal L}(P_{n+1})(x)=x{\mathcal L}(Q_n)(x)-\mathcal L(Q_n^*)(x),
\end{equation}
where
\begin{align*}Q_n^*(x_1,\ldots,x_n)&=\sum_{i_1,\ldots,i_n\in\N}c_{i_1,\ldots,i_n}
x_1x_n\f{\partial^2}{\partial x_1\partial x_n}(x_1^{i_1}\cdots x_n^{i_n})
=x_1x_n\f{\partial^2}{\partial x_1\partial x_n}Q_n(x_1,\ldots,x_n).
\end{align*}
Observe that
\begin{align*}x_1x_n\f{\partial^2}{\partial x_1\partial x_n}Q_n(x_1,\ldots,x_n)
=&\ x_1x_n\f{\partial^2}{\partial x_1\partial x_n}(x_1-x_2)\cdots(x_{n-1}-x_n)
\\=&\ x_1x_n\l((x_2-x_3)\cdots(x_{n-2}-x_{n-1})(-1)\r)
\end{align*}
and thus
$${\mathcal L}(Q_n^*)(x)=-x^2{\mathcal L}(Q_{n-2})(x).$$
Combining this with \eqref{Relation}, we obtain the relation
\begin{equation}
{\mathcal L}(P_{n+1})(x)=x{\mathcal L}(Q_n)(x)+x^2{\mathcal L}(Q_{n-2})(x).
\end{equation}
Note that ${\mathcal L}(Q_{n-2})(x)=x^{(n-3)/2}$ by the induction hypothesis,
and ${\mathcal L}(P_{n+1})(x)=2x^{(n+1)/2}$ by Lemma \ref{Lem2.5}.
Therefore,
$$x{\mathcal L}(Q_n)(x)=2x^{(n+1)/2}-x^2x^{(n-3)/2}=x^{(n+1)/2}$$
and hence ${\mathcal L}(Q_n)(x)=x^{(n-1)/2}$ as desired.

In view of the above, we have completed the induction proof of Lemma \ref{Lem-Q}.

\medskip
\noindent{\it Proof of Theorem} \ref{Th1.3}. Let $k=|A_1|-1=\cdots=|A_n|-1$.
By Lemmas \ref{P*} and \ref{Lem-Q}, we have
\begin{equation}\label{nodd}[x_1^k\cdots x_n^k](x_1-x_2)\cdots(x_{n-1}-x_n)(x_1+\cdots+x_n)^{(k-1)n+1}
=\f{((k-1)n+1)!}{(k!)^n}\times k^{(n-1)/2}.
\end{equation}
Let $h$ be the integer given by the right-hand side of \eqref{nodd}, and let $e$ be the identity of
the field $F$. Clearly, the coefficient of $x_1^k\cdots x_n^k$
in the polynomial
$$(x_1-x_2)\cdots(x_{n-1}-x_n)(x_1+\cdots+x_n)^{(k-1)n+1}\in F[x]$$
is $he$, which is nonzero since $p(F)>(k-1)n+1$. Applying Lemma \ref{Lem1}, we get
$$|L(A_1,\ldots,A_n)|\gs (k-1)n+2=\sum_{i=1}^n|A_i|-2n+2.$$
This concludes our proof. \qed

\bigskip

\section{Proofs of Theorems \ref{Th1} and \ref{Th2}}
\setcounter{lemma}{0}
\setcounter{theorem}{0}
\setcounter{equation}{0}
\setcounter{conjecture}{0}
\setcounter{remark}{0}
\setcounter{corollary}{0}

\medskip
\noindent{\it Proof of Theorem \ref{Th1}}.  The subgroup of $G$ generated by $A$ is a finitely generated torsion-free abelian group.

When $n$ is even, similar to the first paragraph in the proof of Theorem 1.1(ii), without loss of generality we may assume that $G$ is the additive group of the complex field $\C$, hence we obtain
the desired result by applying Theorem 1.2 with $F=\C$ and $A_1=\cdots=A_n=A$.
Similarly, when $n$ is odd we have \eqref{-} by applying Theorem 1.3.

Below we assume that $n$ is odd, and want to prove \eqref{o}.
We may simply assume that $G=\Z^r$ for some positive integer $r$ without any loss of generality.
It is well known that there is a linear ordering $\ls$ on $G=\Z^r$ such that for any $a,b,c\in G$
if $a<b$ then $-b<-a$ and $a+c<b+c$ (cf. \cite{L}). For convenience, we write
$A=\{a_1,a_2,\ldots,a_k\}$ with $a_1<a_2<\ldots< a_k$.

Clearly \eqref{o} holds trivially when $k=2$.
Below we assume that $k\gs3$.

Observe that
\[\begin{aligned}
&\ a_1+a_2+a_1+a_2+\cdots+a_1+a_2+a_3\\
<&\ a_1+a_2+a_1+a_2+\cdots+a_1+a_2+a_4<\cdots\\
<&\ a_1+a_2+a_1+a_2+\cdots+a_1+a_2+a_k
\end{aligned}\]
and
\[\begin{aligned}
&\ a_1+a_i+a_1+a_i+\cdots+a_1+a_i+a_k\\
<&\ a_1+a_{i+1}+a_1+a_i+\cdots+a_1+a_i+a_k<\cdots\\
<&\ a_1+a_{i+1}+a_1+a_{i+1}+\cdots+a_1+a_{i+1}+a_k
\end{aligned}\]
for $2\ls i\ls k-2$. Also,
\[\begin{aligned}
&\ a_1+a_{k-1}+a_1+a_{k-1}+\cdots+a_1+a_{k-1}+a_1+a_{k-1}+a_k\\
<&\ a_1+a_{k}+a_1+a_{k-1}+\cdots+a_1+a_{k-1}+a_1+a_{k-1}+a_k<\cdots\\
<&\ a_1+a_{k}+a_1+a_{k}+\cdots+a_1+a_{k}+a_1+a_{k-1}+a_k,
\end{aligned}\]
and
\[\begin{aligned}
&\ a_i+a_{k}+a_i+a_{k}+\cdots+a_i+a_{k}+a_i+a_{k-1}+a_k\\
<&\ a_{i+1}+a_{k}+a_i+a_{k}+\cdots+a_i+a_{k}+a_i+a_{k-1}+a_k<\cdots\\
<&\ a_{i+1}+a_{k}+a_{i+1}+a_{k}+\cdots+a_{i+1}+a_{k}+a_{i+1}+a_{k-1}+a_k
\end{aligned}\]
for $1\ls i\ls k-3$. Note also that
\[\begin{aligned}
&\ a_{k-2}+a_{k}+a_{k-2}+a_{k}+\cdots+a_{k-2}+a_{k}+a_{k-2}+a_{k-1}+a_k\\
<&\ a_{k-1}+a_{k}+a_{k-2}+a_{k}+\cdots+a_{k-2}+a_{k}+a_{k-2}+a_{k-1}+a_k<\cdots\\
<&\ a_{k-1}+a_{k}+a_{k-1}+a_{k}+\cdots+a_{k-1}+a_{k}+a_{k-2}+a_{k-1}+a_k.
\end{aligned}\]
So we have found
\begin{align*}&1+(k-3)+2\times\l(\f{n-3}2(k-2)+(k-3)\r)
\\=&\ (k-2)n-2=n|A|-2n+(-1)^n(1+\{n\}_2)
\end{align*}
different elements of $n^\circ A$. Therefore \eqref{o} is valid.

In view of the above, we have completed our proof of Theorem \ref{Th1}.
\qed

\medskip
\noindent{\it Proof of Theorem \ref{Th2}}. The subgroup of $G$ generated by $A$
is a finitely generated torsion-free abelian group and hence is isomorphic to $Z^r$ for some positive integer $r$. As there is a linear ordering $\ls$ on $\Z^r$ such that for any $a,b,c\in G$
if $a<b$ then $-b<-a$ and $a+c<b+c$ (cf. \cite{L}), without loss of generality we simply assume that $G=\Z$
and write $A=\{a_1,a_2,\ldots,a_k\}$ with $a_1<a_2<\cdots< a_k$.

(i) If $A$ is an AP with $|A|=k$, then
for $A_0=\{0,\ldots,k-1\}$ we clearly have
$$|n\,\tilde{}A|=|n\,\tilde{}A_0|\ \t{and}\ |n^{\circ}A|=|n^{\circ}{}A_0|.$$
If $n=2m$ with $m\in\Z^+$, then
$$n\,\tilde{}A_0=n^\circ A_0=\{m,m+1,\ldots,m(2k-3)\}$$
with $m=0+1+0+1+\cdots+0+1$ and
$$m(2k-3)=(k-1)+(k-2)+\cdots+(k-1)+(k-2),$$
and hence
$$|n\,\tilde{}A|=|\{m,m+1,\ldots,m(2k-3)\}|=m(2k-3)-(m-1)=2m(k-2)+1=n|A|-2n+1.$$
If $n=2m+1$ with $m\in\Z^+$, then
$$n\,\tilde{}A_0=\{m,m+1,\ldots,m(2k-3)+k-1\}$$
with $m=0+1+0+1+\cdots+0+1+0$ and
$$m(2k-3)=(k-1)+(k-2)+\cdots+(k-1)+(k-2)+(k-1),$$
hence
$$|n\,\tilde{}A|=|\{m,m+1,\ldots,m(2k-3)+k-1\}|
=(2m+1)k-4m=kn-2(n-1)=|A|n-2n+1+\{n\}_2.$$
For $n=2m+1$ with $m\in\Z^+$, we have
$$n^\circ A_0=\{m+2,m+3,\ldots,m(2k-3)+k-3\}$$
with $m=0+1+0+1+\cdots+0+1+2$ and
$$m(2k-3)+k-3=(k-1)+(k-2)+\cdots+(k-1)+(k-2)+(k-1)+(k-3),$$
hence
$$|n^\circ{}A|=|\{m+2,m+3,\ldots,m(2k-3)+k-3\}|
=(2m+1)k-4m-4=|A|n-2n+(-1)^n(1+\{n\}_2).$$
Thus both \eqref{-eq} and \eqref{oeq} hold if $A$ is an AP.

Now we consider the case $n=5$ and $A=\{a_1,a_2,a_3\}$ with $a_1<a_2<a_3$.
Clearly, any element of $5^\circ A$ can be written as $\sum_{k=1}^3n_ka_k$
 with $n_1,n_2,n_3\in\N$ and $n_1+n_2+n_3=5$. Note that $n_k\in\{1,2\}$ for all $k=1,2,3$
 (otherwise we get a contradiction in view of the definition of $5^\circ A$).
 Thus two of $n_1,n_2,n_3$ are $2$ and the remaining one is $1$.
  Hence
  \begin{align*}5^\circ A&=\{a_1+a_2+a_1+a_2+a_3,\ a_1+a_3+a_1+a_2+a_3,\ a_2+a_3+a_1+a_2+a_3\}
  \\&=\{2a_1+2a_2+a_3,\ 2a_1+2a_3+a_2,\ 2a_2+2a_3+a_1\}
  \end{align*}
   contains exactly $3$ elements, and so \eqref{oeq} holds for $k=3$ and $n=5$.

  (ii) Now we write $A=\{a_1,\ldots,a_k\}$ with $a_1<\cdots<a_k$.
  If \eqref{oeq} holds for $n=3$, then $|3^\wedge A|=|3^\circ A|=3k-8$
  and hence $A$ is an AP by Nathanson \cite{N95}.
 Below we
  divide our remaining discussions into four cases.
  \smallskip

  {\it Case} 1. $k\gs4$ and $2\mid n$.

 In this case, both the right-hand sides of \eqref{-eq} and \eqref{oeq}
  are $kn-2n+1$.  As $n\,\tilde{}A\supseteq n^{\circ} A$, if $|n\,\tilde{}A|=kn-2n+1$,
  then $|n^\circ{}A|=kn-2n+1$.

  Now we suppose that $|n^\circ{}A|=kn-2n+1$.
Note that
\[\begin{aligned}
&\ a_1+a_i+a_1+a_i+\cdots+a_1+a_i\\
<&\ a_1+a_{i+1}+a_1+a_i+\cdots+a_1+a_i<\ldots\\
<&\ a_1+a_{i+1}+a_1+a_{i+1}+\cdots+a_1+a_{i+1}
\end{aligned}\]
for all $i\in\{2,3,\ldots,k-1\}$
and
\[\begin{aligned}
&\ a_i+a_k+a_i+a_k+\cdots+a_i+a_k\\
<&\ a_{i+1}+a_k+a_i+a_k+\cdots+a_i+a_k<\ldots\\
<&\ a_{i+1}+a_k+a_{i+1}+a_k+\cdots+a_{i+1}+a_k
\end{aligned}\]
for all $i\in\{1,2,\ldots,k-2\}$.
Therefore we get
$$1+\f{n}2(k-2)\times 2=(k-2)n+1$$ different elements of $n^\circ A$.
They are all the elements of $A$ since $|n^\circ A|=(k-2)n+1$.

For any $i\in\{3,\ldots,k-1\}$, in $n^\circ A$ we have
\[\begin{aligned}
&\ a_1+a_i+a_1+a_i+\cdots+a_1+a_i+a_1+a_{i-1}\\
<&\ a_1+a_{i}+a_1+a_i+\cdots+a_1+a_i+a_1+a_i\\
<&\ a_1+a_{i+1}+a_1+a_{i}+\cdots+a_1+a_i+a_1+a_{i}
\end{aligned}\]
and
\[\begin{aligned}
&\ a_1+a_i+a_1+a_i+\cdots+a_1+a_i+a_1+a_{i-1}\\
<&\ a_1+a_{i+1}+a_1+a_i+\cdots+a_1+a_i+a_1+a_{i-1}\\
<&\ a_1+a_{i+1}+a_1+a_{i}+\cdots+a_1+a_i+a_1+a_{i},
\end{aligned}\]
thus
$$a_1+a_{i}+a_1+a_i+\cdots+a_1+a_i+a_1+a_i=a_1+a_{i+1}+a_1+a_i+\cdots+a_1+a_i+a_1+a_{i-1}$$
and hence  $a_{i+1}-a_i=a_i-a_{i-1}$.
Similarly, for any $i\in\{2,\ldots,k-2\}$, in $n^\circ A$ we have
\[\begin{aligned}
&\ a_i+a_k+a_i+a_k+\cdots+a_i+a_k+a_{i-1}+a_k\\
<&\ a_i+a_k+a_i+a_k+\cdots+a_i+a_k+a_i+a_k\\
<&\ a_{i+1}+a_k+a_i+a_k+\cdots+a_i+a_k+a_i+a_k
\end{aligned}\]
and
\[\begin{aligned}
&\ a_i+a_k+a_i+a_k+\cdots+a_i+a_k+a_{i-1}+a_k\\
<&\ a_{i+1}+a_k+a_i+a_k+\cdots+a_i+a_k+a_{i-1}+a_k\\
<&\ a_{i+1}+a_k+a_i+a_k+\cdots+a_i+a_k+a_i+a_k,
\end{aligned}\]
 and hence $a_{i+1}-a_i=a_i-a_{i-1}$.
 Therefore $a_2-a_1=a_3-a_2=\ldots=a_k-a_{k-1}$, i.e., $A$ is an AP.
\medskip

{\it Case} 2. $k\gs4$ and $2\nmid n$.

Clearly,
\[\begin{aligned}
&\ a_1+a_i+a_1+a_i+\cdots+a_1+a_i+a_1\\
<&\ a_1+a_{i+1}+a_1+a_i+\cdots+a_1+a_i+a_1<\ldots\\
<&\ a_1+a_{i+1}+a_1+a_{i+1}+\cdots+a_1+a_{i+1}+a_1
\end{aligned}\]
for any $i\in\{2,3,\ldots,k-1\}$,
and
\[\begin{aligned}
&\ a_i+a_k+a_i+a_k+\cdots+a_i+a_k+a_i\\
<&\ a_{i+1}+a_k+a_i+a_k+\cdots+a_i+a_k+a_i<\ldots\\
<&\ a_{i+1}+a_k+a_{i+1}+a_k+\cdots+a_{i+1}+a_k+a_{i+1}
\end{aligned}\]
for all $i\in\{1,2,\ldots,k-2\}$. Also,
\begin{align*}
&a_{k-1}+a_k+a_{k-1}+a_k+\cdots+a_{k-1}+a_k+a_{k-1}
\\<&\ a_k+a_{k-1}+a_k+a_{k-1}+\cdots+a_k+a_{k-1}+a_k.
\end{align*}
So we have found
$$1+\f{n+1}2(k-2)+\f{n-1}2(k-2)+1=2+(k-2)n=n|A|-2n+1+\{n\}_2$$ different elements of $n\,\tilde{} A$.
They are all the elements of $n\,\tilde{} A$ if \eqref{-eq} holds.

Now suppose that \eqref{-eq} is valid.
For any $i\in\{2,\ldots,k-2\}$, in $n\,\tilde{} A$ we have
\[\begin{aligned}
&\ a_i+a_k+a_i+a_k+\cdots+a_i+a_k+a_{i-1}\\
<&\ a_i+a_k+a_i+a_k+\cdots+a_i+a_k+a_i\\
<&\ a_{i+1}+a_k+a_i+a_k+\cdots+a_i+a_k+a_i
\end{aligned}\]
and
\[\begin{aligned}
&\ a_i+a_k+a_i+a_k+\cdots+a_i+a_k+a_{i-1}\\
<&\ a_{i+1}+a_k+a_i+a_k+\cdots+a_i+a_k+a_{i-1}\\
<&\ a_{i+1}+a_k+a_i+a_k+\cdots+a_i+a_k+a_i,
\end{aligned}\]
hence
$$a_i+a_k+a_i+a_k+\cdots+a_i+a_k+a_i=a_{i+1}+a_k+a_i+a_k+\cdots+a_i+a_k+a_{i-1}$$
and thus $a_{i+1}-a_i=a_i-a_{i-1}$.
Moreover, in $n\,\tilde{} A$ we have
\[\begin{aligned}
&\ a_{k-1}+a_k+a_{k-1}+a_k+\cdots+a_{k-1}+a_k+a_{k-2}\\
<&\ a_{k-1}+a_k+a_{k-1}+a_k+\cdots+a_{k-1}+a_k+a_{k-1}\\
<&\ a_k+a_{k-1}+a_k+a_{k-1}+\cdots+a_k+a_{k-1}+a_k
\end{aligned}\]
and
\[\begin{aligned}
&\ a_{k-1}+a_k+a_{k-1}+a_k+\cdots+a_{k-1}+a_k+a_{k-2}\\
<&\ a_k+a_{k-1}+a_k+a_{k-1}+\cdots+a_k+a_{k-2}+a_k\\
<&\ a_k+a_{k-1}+a_k+a_{k-1}+\cdots+a_k+a_{k-1}+a_k,
\end{aligned}\]
thus
\begin{align*}&\ a_{k-1}+a_k+a_{k-1}+a_k+\cdots+a_{k-1}+a_k+a_{k-1}
\\=&\ a_k+a_{k-1}+a_k+a_{k-1}+\cdots+a_k+a_{k-2}+a_k
\end{align*}
and hence $a_{k-1}-a_{k-2}=a_k-a_{k-1}$.
\medskip

Assume $n\gs5$. For any $i\in\{3,\ldots,k-1\}$,
in $n^\circ A$,
by the proof of Theorem \ref{Th1} and the equality \eqref{oeq}, there is a unique element of
$n^\circ A$ between
$$a_1+a_i+a_1+a_i+\cdots+a_1+a_i+a_1+a_{i-1}+a_k
\ \t{and}\ a_1+a_{i+1}+a_1+a_{i}+\cdots+a_1+a_i+a_1+a_{i}+a_k.$$
In $n^\circ A$, we clearly have
\[\begin{aligned}
&\ a_1+a_i+a_1+a_i+\cdots+a_1+a_i+a_1+a_{i-1}+a_k\\
<&\ a_1+a_{i}+a_1+a_i+\cdots+a_1+a_i+a_1+a_i+a_k\\
<&\ a_1+a_{i+1}+a_1+a_{i}+\cdots+a_1+a_i+a_1+a_{i}+a_k
\end{aligned}\]
and
\[\begin{aligned}
&\ a_1+a_i+a_1+a_i+\cdots+a_1+a_i+a_1+a_{i-1}+a_k\\
<&\ a_1+a_{i+1}+a_1+a_i+\cdots+a_1+a_i+a_1+a_{i-1}+a_k\\
<&\ a_1+a_{i+1}+a_1+a_{i}+\cdots+a_1+a_i+a_1+a_{i}+a_k,
\end{aligned}\]
Thus $$a_1+a_{i}+a_1+a_i+\cdots+a_1+a_i+a_1+a_i+a_k
=a_1+a_{i+1}+a_1+a_i+\cdots+a_1+a_i+a_1+a_{i-1}+a_k$$
and hence $a_{i+1}-a_i=a_i-a_{i-1}$.

Similarly, for $i\in\{2,\ldots,k-2\}$, by the proof of Theorem \ref{Th1} and the equality \eqref{oeq}, there is a unique element of
$n^\circ A$ between
$$a_i+a_k+a_i+a_k+\cdots+a_i+a_k+a_{i-1}+a_{k-1}+a_k
\ \t{and}\ a_{i+1}+a_k+a_i+a_k+\cdots+a_i+a_k+a_i+a_{k-1}+a_k.$$
In $n^\circ A$ we have
\begin{align*}
&\ a_i+a_k+a_i+a_k+\cdots+a_i+a_k+a_{i-1}+a_{k-1}+a_k\\
<&\ a_i+a_k+a_i+a_k+\cdots+a_i+a_k+a_i+a_{k-1}+a_k\\
<&\ a_{i+1}+a_k+a_i+a_k+\cdots+a_i+a_k+a_i+a_{k-1}+a_k
\end{align*}
and
\[\begin{aligned}
&\ a_i+a_k+a_i+a_k+\cdots+a_i+a_k+a_{i-1}+a_{k-1}+a_k\\
<&\ a_{i+1}+a_k+a_i+a_k+\cdots+a_i+a_k+a_{i-1}+a_{k-1}+a_k\\
<&\ a_{i+1}+a_k+a_i+a_k+\cdots+a_i+a_k+a_i+a_{k-1}+a_k,
\end{aligned}\]
hence
\begin{align*}&\ a_i+a_k+a_i+a_k+\cdots+a_i+a_k+a_i+a_{k-1}+a_k
\\=&\ a_{i+1}+a_k+a_i+a_k+\cdots+a_i+a_k+a_{i-1}+a_{k-1}+a_k
\end{align*}
and thus
 $a_{i+1}-a_i=a_i-a_{i-1}$.
Therefore $a_2-a_1=a_3-a_2=\ldots=a_k-a_{k-1}$, i.e., $A$ is an AP.
\medskip

{\it Case} 3. $k=3$ and $2\mid n$.

If \eqref{-eq} holds then so \eqref{oeq}.
Now suppose that \eqref{oeq} holds.
In $n^\circ A$, we have
\[\begin{aligned}
&\ a_1+a_2+a_1+a_2+\cdots+a_1+a_2\\
<&\ a_1+a_{3}+a_1+a_2+\cdots+a_1+a_2
< \cdots\\
<&\ a_1+a_{3}+a_1+a_{3}+\cdots+a_1+a_{3}
\\
<&\ a_{2}+a_3+a_2+a_3+\cdots+a_2+a_3
< \cdots\\
<&\ a_{2}+a_3+a_{2}+a_3+\cdots+a_{2}+a_3,
\end{aligned}\]
they give all the $n+1$ elements of $n^\circ A$.
As
\[\begin{aligned}
&\ a_1+a_3+a_1+a_3+\cdots+a_1+a_3+a_1+a_2\\
<&\ a_1+a_3+a_1+a_3+\cdots+a_1+a_3+a_1+a_3\\
<&\ a_2+a_3+a_1+a_3+\cdots+a_1+a_3+a_1+a_3
\end{aligned}\]
and
\[\begin{aligned}
&\ a_1+a_3+a_1+a_3+\cdots+a_1+a_3+a_1+a_2\\
<&\ a_2+a_3+a_1+a_3+\cdots+a_1+a_3+a_1+a_2\\
<&\ a_2+a_3+a_1+a_3+\cdots+a_1+a_3+a_1+a_3,
\end{aligned}\]
we must have
$$a_1+a_3+a_1+a_3+\cdots+a_1+a_3+a_1+a_3=a_2+a_3+a_1+a_3+\cdots+a_1+a_3+a_1+a_2,$$
and hence $a_3-a_2=a_2-a_1$.
\medskip

{\it Case} 4. $k=3$ and $2\nmid n$.

In $n\,\tilde{} A$, we have
\begin{align*}
&\ a_1+a_2+a_1+a_2+\cdots+a_1+a_2+a_1\\
<&\ a_1+a_{3}+a_1+a_2+\cdots+a_1+a_2+a_1<\cdots\\
<&\ a_1+a_{3}+a_1+a_{3}+\cdots+a_1+a_{3}+a_1
\\
<&\ a_{2}+a_3+a_1+a_3+\cdots+a_1+a_3+a_1<\cdots\\
<&\ a_{2}+a_3+a_{2}+a_3+\cdots+a_{2}+a_3+a_{2}
\\<&\ a_3+a_{2}+a_3+a_{2}+\cdots+a_3+a_{2}+a_3.
\end{align*}
Suppose that \eqref{-eq} holds. Then the above gives a list of all the $n+2$ elements of $n\,\tilde{}A$.

When $n=3$, in $n\,\tilde{} A$ we have
\[a_2+a_3+a_1<a_2+a_3+a_2<a_3+a_2+a_3
\]
and
\[a_2+a_3+a_1<a_3+a_1+a_3<a_3+a_2+a_3,
\]
hence $a_2+a_3+a_2=a_3+a_1+a_3$  and thus
$a_2-a_1=a_3-a_2$.
When $n\gs5$, as
\[\begin{aligned}
&\,a_1+a_3+a_1+a_3+\cdots+a_1+a_3+a_1+a_2+a_1\\
<&\,a_1+a_3+a_1+a_3+\cdots+a_1+a_3+a_1+a_3+a_1\\
<&\,a_2+a_3+a_1+a_3+\cdots+a_1+a_3+a_1+a_3+a_1
\end{aligned}\]
and
\[\begin{aligned}
&\,a_1+a_3+a_1+a_3+\cdots+a_1+a_3+a_1+a_2+a_1\\
<&\,a_2+a_3+a_1+a_3+\cdots+a_1+a_3+a_1+a_2+a_1\\
<&\,a_2+a_3+a_1+a_3+\cdots+a_1+a_3+a_1+a_3+a_1,
\end{aligned}\]
we must have
$$a_1+a_3+a_1+a_3+\cdots+a_1+a_3+a_1+a_3+a_1=a_2+a_3+a_1+a_3+\cdots+a_1+a_3+a_1+a_2+a_1$$
and hence
$a_2-a_1=a_3-a_2$.

Assume $n\gs7$ and suppose that \eqref{oeq} holds.
By the proof of Theorem \ref{Th1} and the equality \eqref{oeq}, there is a unique element of
$n^\circ A$ between
$$a_1+a_3+a_1+a_3+\cdots +a_1+a_3+a_1+a_2+a_1+a_2+a_3$$
and $$a_2+a_3+a_1+a_3+\cdots +a_1+a_3+a_1+a_3+a_1+a_2+a_3.$$
In $n^\circ A$, we clearly have
\[\begin{aligned}
&\ a_1+a_3+a_1+a_3+\cdots +a_1+a_3+a_1+a_2+a_1+a_2+a_3\\
<&\ a_1+a_3+a_1+a_3+\cdots +a_1+a_3+a_1+a_3+a_1+a_2+a_3\\
<&\ a_2+a_3+a_1+a_3+\cdots +a_1+a_3+a_1+a_3+a_1+a_2+a_3
\end{aligned}\]
and
\[\begin{aligned}
&\ a_1+a_3+a_1+a_3+\cdots +a_1+a_3+a_1+a_2+a_1+a_2+a_3\\
<&\ a_2+a_3+a_1+a_3+\cdots +a_1+a_3+a_1+a_2+a_1+a_2+a_3\\
<&\ a_2+a_3+a_1+a_3+\cdots +a_1+a_3+a_1+a_3+a_1+a_2+a_3.
\end{aligned}\]
Thus
\begin{align*}&a_1+a_3+a_1+a_3+\cdots +a_1+a_3+a_1+a_3+a_1+a_2+a_3
\\=&\ a_2+a_3+a_1+a_3+\cdots +a_1+a_3+a_1+a_2+a_1+a_2+a_3,
\end{align*}
and hence $a_2-a_1=a_3-a_2$ as desired.

In view of the above, we have finished our proof of Theorem \ref{Th2}.
\qed

\end{document}